\documentclass[scale=1.1]{article}
\usepackage[margin=1in]{geometry}
\usepackage{url}
\usepackage{plain}
\usepackage[T1]{fontenc}\usepackage{lmodern}

\usepackage{float}
\usepackage{tikz}
\usepackage{subcaption}
\usetikzlibrary{tikzmark,decorations.pathreplacing}
\usetikzlibrary{arrows,matrix}
\usetikzlibrary{positioning}
\usetikzlibrary{decorations}
\usetikzlibrary{calc}

\begin{document}
\begin{plain}

%begin macros
\def\L{{\cal L}}
\baselineskip=14pt
\parskip=10pt
\def\halmos{\hbox{\vrule height0.15cm width0.01cm\vbox{\hrule height
  0.01cm width0.2cm \vskip0.15cm \hrule height 0.01cm width0.2cm}\vrule
  height0.15cm width 0.01cm}}
\font\eightrm=cmr8 \font\sixrm=cmr6
\font\eighttt=cmtt8
%\magnification=\magstephalf
\def\P{{\cal P}}
\def\Q{{\cal Q}}
\def\S{{\cal S}}
\def\1{{\overline{1}}}
\def\2{{\overline{2}}}
\parindent=0pt
\overfullrule=0in
\def\Tilde{\char126\relax}
\def\frac#1#2{{#1 \over #2}}
%\headline={\rm  \ifodd\pageno  \RightHead  \else  \LeftHead  \fi}
%\def\RightHead{\centerline{
%Title
%}}
%\def\LeftHead{ \centerline{Doron Zeilberger}}
%end macros
\centerline
{\bf
In How Many Ways can a Rectangle be Rectangled?
}
\bigskip
\centerline
{\it Pablo BLANCO, Robert DOUGHERTY-BLISS, Natalya TER-SAAKOV, and Doron ZEILBERGER}

\bigskip

{\bf Abstract}: There are $2^{n-1}$  ways to tile a $1 \times n$ rectangle with rectangular tiles (of any length, of course they all must have width $1$), but
in how many ways can you tile  a $100 \times 100$  checkerboard with such tiles? Neither humankind, nor computer-kind, will (most probably) ever know the exact number.
But it is possible to compute these numbers for $m \times n$ rectangular grids, if $m$ is not too big, while $n$ can be as big as one wishes.
This was initially done in 1988 by David  Klarner and  Spyros Magliveras, and
beautifully extended, around 2006, by, at-the-time, first-year LSU undergraduate  Joshua Smith, in collaboration with his faculty mentor, Helena Verrill.
Here we extend this to weighted-counting, also keeping track of the number of tiles (that ranges from $1$ to $mn$), and
the number of participating grid-edges (that range from $2m+2n$ to $2mn+m+n$).
This quickly leads to statistical analyses (mean, variance, and higher moments) of  these quantities.
While we admire the clever approaches of Klarner-Magliveras and Smith-Verrill, we use two alternative approaches to the original problem,
that are more amenable for deriving these generalizations. At the same time, we illustrate the power and beauty of experimental-yet-rigorous enumerative combinatorics.

{\bf Preface: How it all started}

A few weeks ago, the New York Times Sunday magazine started publishing a puzzle, composed by Prasanna Seshadri, that they call {\it Recto}.
You are given a rectangular grid, say $m \times n$  (in their case it is always a square, i.e. $m=n$), where some of the boxes contain positive integers larger than one.
In their own words, you have to  do as follows:

{\it ``Divide the grid into rectangles (including squares) so that each region contains one number. The number will represent the sum of the length and height of that region.
Regions may not overlap and all cells must be used.''}

See here for a  sample of solved puzzles:

{\tt \url{https://sites.math.rutgers.edu/~zeilberg/mamarim/mamarimhtml/RectoSample.pdf}} \quad .

So the solution is always a certain tiling of the $m \times n$   checkerboard by rectangular tiles of {\it any} dimensions. 

Being enumerators, the natural question that came to our mind was: {\it How many such tilings are there?}
Let's call that number $a(m,n)$.
We are almost sure, even with AI getting smarter and smarter, that neither humankind nor machinekind will ever know the {\it exact} value of $a(100,100)$,
but {\it we'll do what we can.}  Our goal was to find explicit expressions, in the variable $x$, for the generating functions $\sum_{n=0}^{\infty} a(m,n)\,x^n$, for as large $m$ as the
computer would allow. So we wrote a Maple package, {\it RectTile.txt}, available from

{\tt \url{https://sites.math.rutgers.edu/~zeilberg/tokhniot/RectTile.txt}} \quad .

{\bf First Step: Crank out a few values}

Using {\it Dynamical programming}, it is not too hard to find these numbers for moderate-size rectangular grids.

Suppose that you want to find that number for an $m \times n$ grid, for a fixed $m$ but general $n$. We agreed to denote that number by $a(m,n)$.
The top-right box must be part of {\it some} rectangular tile
of dimension $a \times b$, say, for some $1\leq a \leq m$, and some $1 \leq b \leq n$. Removing it leaves us with a smaller board, but no longer rectangular.
So we are {\bf forced} to consider more general boards. By always  looking at the rightmost-topmost still surviving box, and
all the possible ways of removing a rectangular tile containing it, and keep going, our intermediate boards
can be described by vectors of non-negative integers
$$
[c_1, \dots, c_m] \quad,
$$
where, for each row $i$, $1\leq i \leq m$,  $c_i$ (that may be zero) is the number of surviving boxes in the $i$-th row, in the
(necessarily) left-justified current board.

We call these {\it configurations}. There is a natural notion of {\it children}, that consists
of  the set of all boards obtained from it by removing a rectangular tile that contains the top-right box and
that is fully contained in the configuration.

The {\bf base case} is $[0,0, \dots, 0]$, the empty board, whose number of tilings is $1$ (do nothing).
If $F(L)$ is the number of tilings of the configurations $L$, then of course
$$
F(L) \, = \, \sum_{L' \in Children(L)} F(L') \quad.
$$

At the {\it end of the day}, $a(m,n)$, our desired number, equals
$$
F([n,n, \dots , n]) \quad ,
$$
where $n$ is repeated $m$ times.

{\eightrm

This is implemented in procedure  {\tt TIr(L)} in our Maple package {\tt RectTile.txt}. For example to get the number of  such tilings for a $6 \times 6$ grid, type:

{\tt TIr([6,6,6,6,6,6]);} \quad,

immediately getting that the exact number is: $535236230270$ \quad.
}

But, so far this is {\it numerics}, aka {\it number crunching}. It would be too much to hope for a `formula', or even a polynomial time algorithm in $m+n$, for general $m$ and $n$,
but how about an explicit formula, or just-as-good, an explicit generating function, for the (singly) infinite sequence $\{a(m,n)\}_{n=0}^{\infty}$, when  $m$ is small?
And if we are in luck, when $m$ is not so small.

{\bf Second Step:  Go to the OEIS!}

Whenever enumerators  encounter an  integer sequence new-to-them, they immediately go to
Neil Sloane's  {\it On-Line Encyclopedia of Integer Sequences} (OEIS) [Sl], the most useful (and fascinating!)  mathematical database.
This way they find out whether their sequence is novel. That's exactly what we did with this problem.

For  $m=3$, using the {\it numeric} procedure {\tt TIr}, we typed, in our Maple package,

{\tt seq(TIr([n,n,n]),n=1..11);} \quad,

and immediately got:
$$
4, 34, 322, 3164, 31484, 314662, 3149674, 31544384, 315981452, 3165414034, 31710994234
$$

Then we checked the OEIS for the above numbers and, to our initial dismay, it turned out that we were scooped!

Sure enough this sequence has been there since April 22, 2012 (and contributed by our mathematical son (for DZ) and mathematical brother (for PB and RDB), the
brilliant Mathew C. Russell.  It is OEIS sequence  {\tt https://oeis.org/A208215} .

Moving right along we found that (all entered by Alois P. Heinz in Dec. 2012)

$\bullet$  for $m=4$ we got {\tt https://oeis.org/A220297} \quad ;

$\bullet$  for $m=5$ we got {\tt https://oeis.org/A220298} \quad ;

$\bullet$  for $m=6$ we got  {\tt https://oeis.org/A220299} \quad ;

$\bullet$  for $m=7$ we got  {\tt https://oeis.org/A220300} \quad ;

$\bullet$  for $m=8$ we got  {\tt https://oeis.org/A220301} \quad ;

$\bullet$  for $m=9$ we got  {\tt https://oeis.org/A220302} \quad ;

$\bullet$  for $m=10$ we got  {\tt https://oeis.org/A220303} \quad .

The entries for $3 \leq m \leq 6$ all contained the generating functions of the sequences, but for $7 \leq m \leq 10$ they were absent (viewed June 1, 2026).

Finally  {\tt https://oeis.org/A116694} combines them all into a `triangle'. This was contributed by Helena Verrill on Feb. 13, 2006. 

So indeed, as far, as {\it just} counting, we have been scooped. Also thanks to the OEIS, we found a real gem [SV], extending  a paper by
David Klarner and  Spyros Magliveras [KS]. The OEIS did a great service of making [SV] public by
producing a {\it cached} version. As far as we can tell, it is neither published in a peer-reviewed journal, nor can it be found on arxiv.org \quad .

According to a footnote of [SV], this paper was written while the first author, Joshua Smith, was a first-year undergraduate(!), presumably advised by
then LSU professor, Helena Verrill, now in Warwick. We are very impressed. 

They proved:

{\bf Theorem} (J. Smith and H. Verrill, [SV], Theorem 1)
$$
a(m,n) \, = \, {\bf 1} \cdot (M_m)^{n-1} \cdot {\bf 1}^T \quad,
$$
where ${\bf 1}=(1,\dots, 1)\in Z^{2^{m-1}}$, and $M_m$ is a $2^{m-1} \times 2^{m-1}$ matrix, defined recursively as follows
$$
M_1=(2) \quad, \quad B_1=(1) \quad, \quad
M_{m+1}= \left ( \matrix{ M_m & B_m \cr
                          B_m & 2M_m} \right ) \quad, \quad
B_{m+1}= \left ( \matrix{ B_m & B_m \cr
                          B_m & M_m} \right )  \quad .
$$

So we immediately have:

{\bf Corollary}: For a fixed $m$, the generating function of the sequence $\{a(m,n) \}_{n=0}^{\infty}$, let's call it $F_m(x)$, i.e.
$$
F_m(x) \,:= \,  \sum_{n=0}^{\infty} a(m,n) x^n \quad,
$$
is given by the formula
$$
F_m(x) \, = 1+ \, x {\bf 1} (I \,- \, x M_m)^{-1} {\bf 1}^T \quad .
$$

This is implemented in procedure  {\tt FmxSV(m,x)} in our Maple package. This produced the following explicit expressions, confirming the ones in the OEIS for $2 \leq m \leq 6$, and
stating (possibly) for the first time, the rational function expressions for $7 \leq m \leq 9$.

$$
F_2(x) \, = \, \frac{3 x^{2}-4 x +1}{7 x^{2}-6 x +1} \quad ,
$$
$$
F_3(x) \, = \, \frac{19 x^{3}-29 x^{2}+11 x -1}{51 x^{3}-55 x^{2}+15 x -1}  \quad ,
$$
$$
F_4(x) \, = \, \frac{3832 x^{6}-8492 x^{5}+6722 x^{4}-2468 x^{3}+441 x^{2}-36 x +1}{11680 x^{6}-20980 x^{5}+13840 x^{4}-4280 x^{3}+645 x^{2}-44 x +1} \quad ,
$$
$$
F_5(x) \, = \, \frac{N_5(x)}{D_5(x)} \quad,
$$
%\vfill\eject
where
$$
N_5(x):= \, 39672144 x^{10}-110891556 x^{9}+124284414 x^{8}-74544838 x^{7}+26669637 x^{6}-5961522 x^{5} 
$$
$$
+ 841659 x^{4}-73608 x^{3}+3769 x^{2}-100 x +1 \quad,
$$
$$
D_5(x):=\,  135762480 x^{10}-326041524 x^{9}+320708934 x^{8}-170972730 x^{7}+54776249 x^{6}
$$
$$
-11002298 x^{5}+1395665 x^{4}-109292 x^{3}+4975 x^{2}-116 x +1 \quad .
$$

The expressions for $6 \leq m \leq 9$ get more and more complicated. They can be viewed in this output file:

{\tt \url{https://sites.math.rutgers.edu/~zeilberg/tokhniot/oRectTile1sv.txt}} \quad .

{\bf Getting  Generating Functions by Guessing the Answer}

The Smith-Verrill method is very clever and elegant, but what if we are not so smart? Then we can easily get the first few members of $F_m(x)$, by plain {\it guessing}.
Using the above-mentioned all-purpose procedure {\tt TIr(L)}, we ask our computer to crank out sufficiently many terms and then fit the
data to the $C-finite$ ansatz [Z], in other words use {\it undetermined coefficients} to find a rational function whose Taylor coefficients match our sequence.
We actually implemented this naive approach in procedure {\tt FmxG(m,x)}. We got the same output, but it took longer.
Yet the agreement was reassuring.

{\bf Getting  Generating Functions by Guessing the Grammar}

To every tiling of an $m \times n$ rectangular grid, we can associate the set of edges that belong to the  boundaries of the participating tiles. So each such tiling corresponds
to a certain subset of the set of $m(n+1)+n(m+1)$ edges of the $m \times n$ checkerboard. Of course, most such subsets do not correspond to any tiling.
Let's label the $(m+1)(n+1)$ {\it vertices}  of the grid in matrix notation.
$$
\{[i,j] \, \Bigl{\vert}  \, 0 \leq i \leq m , \quad  0 \leq j \leq n \} \quad .
$$
The set of edges in the $m \times n$ grid is
$$
\{ \{[i,j], [i+1,j] \} \, \Bigl{\vert} \, 0 \leq i \leq m-1, \quad 0 \leq j \leq n \} \bigcup
\{ \{[i,j], [i,j+1]\} \, \Bigl{\vert}\, 0 \leq i \leq m, \quad 0 \leq j \leq n-1 \} \quad .
$$

So every tiling of the $m \times n$ checkerboard is uniquely determined by the set of participating edges. This can be naturally viewed as
a word of length $n$ with a certain alphabet that consists of {\it some} of the members of the set of edges, let's call it $S(m)$:
$$
S(m):=\{ \{ [i,0],[i+1,0] \} \Bigl{\vert}, 0 \leq i \leq m-1 \} \bigcup
\{ \{ [i,0],[i,1] \} \Bigl{\vert}, 0 \leq i \leq m \}  \quad .
$$

Of course most of the $2^{2m+1}$ subsets never show up. So we have a {\bf language}, and it is reasonable to conjecture that it is a {\it regular language} (aka {\it finite automata}).
In other words, there is a certain set of {\it starting letters}, and then a larger set, (superset of the former), of {\it middle letters}, and each `letter'
may be followed by a certain subset of the letters. Finally the {\bf final} letter is always the `letter' consisting of all the edges of the rightmost border.

For $1 \leq j \leq n$, the $j$-th letter is obtained by intersecting the tiling with the `cross-section' $[0,m] \times [j-1,j)$, and subtracting $j$ from the second coordinates, making it
a subset of $S(m)$ defined above.

Later on we will show how to derive the alphabet and the (type three) grammar using {\it human ingenuity}, but it is also fun to emulate {\it ChatGPT} and
completely {\it without thinking}, conjecture an `alphabet' and a (regular) grammar.

For a fixed $m$, we generate the {\it corpus} of all such $n$-letter words, but $n=3$ suffices (after we converted the tilings into such words), and
let the computer discover {\it all by itself}, {\bf without any `thinking'}:

$\bullet$ The set of starting letters \quad ;

$\bullet$ The set of all letters (that happened to contain the set of starting letters) \quad ;

$\bullet$ For each letter in the above alphabet, which letters can follow?

{\it A priori} it not guaranteed that such a grammar exists (it can be proved, using {\it thinking} see below), but who cares?
This hypothesis can be tested empirically.

The natural set-up is a {\bf directed graph} where the vertices are labeled by the members of our (empirically derived) alphabet, and
there is an edge from the vertex corresponding to the letter $L_1$ to the vertex corresponding to the letter $L_2$ if and only if,
in our language, $L_2$ can come right after $L_1$. We also add two `artificial edges' {\tt START} and {\tt END}, and
put edges between {\tt START} and the starting letters, and between all letters and {\tt END}. Our beloved
computer can construct this {\bf directed graph} fully automatically. Now the set of tilings of the $m \times n$ grid graph is in bijection
with {\it walks} of length $n$ in the directed graph.

Recall that for a directed graph on $N+2$ vertices with a starting vertex labeled $1$ and terminal vertex labeled $N+2$,
if $A$ is the {\it adjacency matrix}, then if $w(n)$ is the number of walks of length $n$, from START($1$) to END ($N+2$),in that graph, then
$$
\sum_{n=0}^{\infty} \, w(n)\,x^n= \left( (I-xA)^{-1} \right )_{1,N+2} \quad .
$$

{\eightrm In our Maple package, {\tt Corpus(m,n)}, generates all the `words' (with our convention). 
In order to get the grammar fully empirically, we have a procedure called {\tt GrammarE(m)}. For example to empirically guess
the alphabet and grammar of rectangular tilings of the $3 \times n$ checkerboard, type {\tt GrammarE(3);}.

Procedure {\tt DiG(m)} converts it to a directed graph, and then using the general procedure {\tt WalkGF(G,x)} we
get the desired generating function.
}

{\bf The Deductive Approach to Deriving the Grammar}

Being broad-minded experimental mathematicians, who care little for that straitjacket called {\it rigor} (that held mathematics back for so long), we really like the
{\it empirical} approach for discovering the alphabet and grammar. Alas, as $m$ grows larger, the corpus
gets larger and larger, so we have to put on our thinking cap and teach the computer how, for a given $m$, discover the alphabet and the grammar
{\it deductively}, rather than {\it inductively}.

Yet, the corpus did help to study the structure of the language.

Using the above convention that a letter is a subset of the set of edges of $[0,m]\times [0,1)$, i.e. of $S(m)$, defined above,
it is easy to see that a starting letter must contain all the vertical edges of the leftmost boundary:
$$
\{\{ [i,0],[i+1,0] \} \Bigl{\vert} 0 \leq i \leq m-1 \}  \quad,
$$
and also the top and bottom horizontal edges
$$
\{[0,0],[0,1]\}, \{[m,0],[m,1]\} \quad,
$$
and {\it any} subset (including the empty set) of the remaining horizontal edges of the leftmost cross-section of the grid
$$
\{ \{ [i,0],[i,1] \} \Bigl{\vert}, 1 \leq i \leq m-1 \} \quad .
$$

Note that there are exactly $2^{m-1}$ {\it starting letters}.

The next thing to {\it teach the computer} is:

{\it Given a letter $L$, what letters can follow it?}

Every letter can be broken into a pair $L=[V,H]$, where $V$ is the
subset of vertical edges and $H$ is the subset of horizontal edges. It is easy to see that the set of possible followers of $L$
{\bf only} depend  on $H$. Let $L'=[V',H']$ be a legal follower of $L=[V,H]$. It is readily seen that any consecutive run of
vertical edges of $V'$ must start and end at one of the endpoints of the edges of $H$.
We ask our beloved computer to generate all possible such $V'$.
Now we look at each section of the “wall” created by V' , and except for the bottom and top (mandatory) edges, you can put any subset that “sticks” out of the section. In the gaps between sections, H' should have an edge if and only if H did.

{\eightrm This is implemented in procedure {\tt Followers(L)} \quad .}

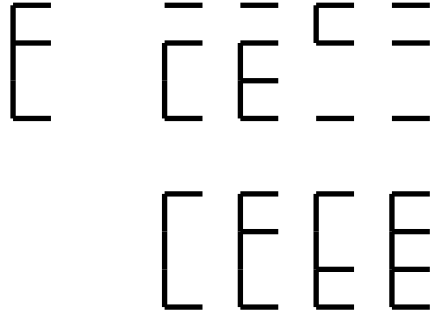
\begin{figure}[H]
    \centerline{
    \begin{tikzpicture}[scale=0.5]

%\tikzstyle{knode}=[rectangle,fill,draw=black,thick,inner sep=1pt] % to make white vertices, erase fill

% grid graph
\pgfmathsetmacro{\n}{4} % set number of vertices

%%% first letter
% horizontal lines
\foreach \y in {1,3,4}
{
 \draw[line width=2pt] (0,\y) -- (1,\y);
}
% vertical lines
\foreach \y in {1,2,3}
{
 \draw[line width=2pt] (0,\y) -- (0,{\y+1});
}

%%% second letter
\pgfmathsetmacro{\L}{2}
\pgfmathsetmacro{\offs}{-1}
% horizontal lines
\foreach \y in {1,4}
{
 \draw[line width=2pt] ({0+2*\L},{\y+5*\offs}) -- ({1+2*\L},{\y+5*\offs});
}
% vertical lines
\foreach \y in {1,2,3}
{
 \draw[line width=2pt] ({0+2*\L},{\y+5*\offs}) -- ({0+2*\L},{\y+1+5*\offs});
}

%%% third letter
\pgfmathsetmacro{\L}{3}
% horizontal lines
\foreach \y in {1,3,4}
{
 \draw[line width=2pt] ({0+2*\L},{\y+5*\offs}) -- ({1+2*\L},{\y+5*\offs});
}
% vertical lines
\foreach \y in {1,2,3}
{
 \draw[line width=2pt] ({0+2*\L},{\y+5*\offs}) -- ({0+2*\L},{\y+1+5*\offs});
}

%%% fourth letter
\pgfmathsetmacro{\L}{4}

% horizontal lines
\foreach \y in {1,2,4}
{
 \draw[line width=2pt] ({0+2*\L},{\y+5*\offs}) -- ({1+2*\L},{\y+5*\offs});
}
% vertical lines
\foreach \y in {1,2,3}
{
 \draw[line width=2pt] ({0+2*\L},{\y+5*\offs}) -- ({0+2*\L},{\y+1+5*\offs});
}

%%% fifth letter
\pgfmathsetmacro{\L}{5}
% horizontal lines
\foreach \y in {1,2,3,4}
{
 \draw[line width=2pt] ({0+2*\L},{\y+5*\offs}) -- ({1+2*\L},{\y+5*\offs});
}
% vertical lines
\foreach \y in {1,2,3}
{
 \draw[line width=2pt] ({0+2*\L},{\y+5*\offs}) -- ({0+2*\L},{\y+1+5*\offs});
}

%%% sixth letter
\pgfmathsetmacro{\L}{2}
% horizontal lines
\foreach \y in {1,3,4}
{
 \draw[line width=2pt] ({0+2*\L},\y) -- ({1+2*\L},\y);
}
% vertical lines
\foreach \y in {1,2}
{
 \draw[line width=2pt] ({0+2*\L},\y) -- ({0+2*\L},{\y+1});
}

%%% seventh letter
\pgfmathsetmacro{\L}{3}
% horizontal lines
\foreach \y in {1,2,3,4}
{
 \draw[line width=2pt] ({0+2*\L},\y) -- ({1+2*\L},\y);
}
% vertical lines
\foreach \y in {1,2}
{
 \draw[line width=2pt] ({0+2*\L},\y) -- ({0+2*\L},{\y+1});
}

%%% eighth letter
\pgfmathsetmacro{\L}{4}
% horizontal lines
\foreach \y in {1,3,4}
{
 \draw[line width=2pt] ({0+2*\L},\y) -- ({1+2*\L},\y);
}
% vertical lines
\foreach \y in {3}
{
 \draw[line width=2pt] ({0+2*\L},\y) -- ({0+2*\L},{\y+1});
}

%%% ninth letter
\pgfmathsetmacro{\L}{5}
\pgfmathsetmacro{\offs}{0}
% horizontal lines
\foreach \y in {1,3,4}
{
 \draw[line width=2pt] ({0+2*\L},{\y+5*\offs}) -- ({1+2*\L},{\y+5*\offs});
}

%%%

\end{tikzpicture}}
    \caption{To the left, one example of a starting letter of a $3\times n$ grid ($n\ge 3$). To the right, the eight potential letters in the grammar which follow that starting letter.}
    \label{fig:grammar example}
\end{figure}

Equipped with {\tt Followers}, we can dynamically create the alphabet, start with the starting letters, and keep applying {\tt Followers} until you
don't encounter any new letters. It turns out that all the letters show up right away. This explains why in the empirical approach,
{\tt Corpus(m,3)} sufficed. Once you have all the letters ({\eightrm implemented in procedure {\tt AlefBet(m)}}),
the computer automatically creates the directed graph, and we proceed as before.

To get $F_m(x)$ using {\it our} approach, type {\tt Fmx(m,x);}. We must admit that the Smith-Verrill approach seems faster, but with our grammar
we can do {\it weighted counting}.

{\bf Weighted Counting}

Our `grammatical' approach is useful for {\it weighted counting}.  let $TILINGS(m,n)$ be the set of tilings by rectangles of the grid-graph $[0,m] \times [0,n]$.
Let
$$
a_T(m,n;t):= \sum_{T \in TILINGS(m,n)} t^{NumberOfTiles(T)} \quad.
$$
$$
a_E(m,n;w):= \sum_{T \in TILINGS(m,n)} w^{NumberOfGridEdges(T)} \quad.
$$
$$
a_{EE}(m,n;w_1,w_2):= \sum_{T \in TILINGS(m,n)} w_1^{NumberOfHorizomalGridEdges(T)} \, w_2^{NumberOfVerticallGridEdges(T)} \quad.
$$

Our approach lets us computer the {\it weighted generating functions}
$$
F^{(T)}_m(x,t):= \sum_{n=0}^{\infty} a_T(m,n;t) \, x^n \quad ,
$$
$$
F^{(E)}_m(x,w):= \sum_{n=0}^{\infty} a_E(m,n;w) \, x^n \quad \,
$$
$$
F^{(EE)}_m(x,w_1,w_2):= \sum_{n=0}^{\infty} a_{EE}(m,n;w_1,w_2) \, x^n \quad .
$$

This is just a minor tweak. Instead of directed graphs, we have {\it weighted} directed graphs where the edges carry weights that keep track
of the desired quantities.

{\eightrm This is implemented in procedures {\tt DiGt(m)}, {\tt DiGw(m)}, and {\tt DiGw1w2(m)}} \quad.

So now instead of the {\it numeric} {\it adjacency matrix}, whose entries are  $0$s and $1$s, we have a matrix of {\it monomials} (and $0$s),
and the formula
$$
\sum_{n=0} \, w(n)\,x^n= \left( (I-xA)^{-1} \right )_{1,N+2} \quad ,
$$
is still applicable, but now $A$ is the appropriate matrix of monomials.

{\eightrm The corresponding procedures are {\tt Fmxt(n,x,t),Fmxw(n,x,w),Fmxw1w2(n,x,w1,w2)} } \quad .

{\bf Sample Output}

$\bullet$ To see $F^{(T)}_m(x,t)$ for $1 \leq m \leq 6$ look here:

 \url{https://sites.math.rutgers.edu/~zeilberg/tokhniot/oRectTile4.txt} \quad .

$\bullet$ To see $F^{(E)}_m(x,w)$ for $1 \leq m \leq 6$ look here:

{\tt \url{https://sites.math.rutgers.edu/~zeilberg/tokhniot/oRectTile3.txt}} \quad .

$\bullet$ To see $F^{(EE)}_m(x,w_1,w_2)$ for $1 \leq m \leq 5$ look here:

{\tt \url{https://sites.math.rutgers.edu/~zeilberg/tokhniot/oRectTile5.txt}} \quad .

For your convenience, these are also available in the package itself, by typing for $1 \leq m \leq 6$:

{\eightrm {\tt FmxtPC(n,x,t),Fmxw(n,x,w)PC,Fmxw1w2PC(n,x,w1,w2)} } \quad .

Also, the pre-computed version of {\tt Fmx(m,x)},  {\tt FmxPC(m,x)}, for {\tt m} between $1$ and $9$ outputs these important generating functions right away.

{\bf Asymptotic Statistical  Analysis}

Once we have explicit weighted  generating functions, we can fully automatically, using Maple's powerful symbolic computation, find out
the asymptotic {\it average number of tiles},  {\it average  number of edges}, and {\it correlation of the joint statistics [the number of horizontal edges,the number of vertical edges]}.

{\eightrm These are  handled by procedures {\tt PaperAsyTiles}, {\tt PaperAsyEdges}, and {\tt PaperAsyCor}, respectively.}

See the output file:

{\tt \url{https://sites.math.rutgers.edu/~zeilberg/tokhniot/oRectTile6.txt}} \quad .

{\bf Random Generation of Tilings}

The Dynamical programming procedure described at the beginning of this paper, {\tt TiR(L)}, enables us to use
Herb Wilf's [W] methodology to generate, uniformaly-at-random, a `typical' tiling, expressed
in the original data structure. For example to get one of the $535236230270$ tilings of
the $6 \times 6$ grid type:

{\tt RandRT([6,6,6,6,6,6]);} \quad.

Let's test it for the $3 \times 3$ case, where there are $322$ such tilings. Typing, in  {\tt RectTile.txt},

{\tt f:=add(X[RandRT([3,3,3])],i=1..32200):} \quad ,

followed by

{\tt add(y**op(1,op(i,f)),i=1..nops(f));} \quad ,

gave us (of course, every time you would get something different)
$$
y^{128}+2 y^{122}+2 y^{121}+3 y^{120}+y^{119}+3 y^{118}+5 y^{117}+4 y^{116}+7 y^{115}+4 y^{114}+8 y^{113}
$$
$$
+3 y^{112}+6 y^{111}+7 y^{110}+7 y^{109}+11 y^{108}+11 y^{107}+11 y^{106}+11 y^{105}+15 y^{104}+6 y^{103}
$$
$$
+13 y^{102}+16 y^{101}+16 y^{100}+11 y^{99}+12 y^{98}+10 y^{97}+8 y^{96}+7 y^{95}+13 y^{94}+11 y^{93}+11 y^{92}
$$
$$
+14 y^{91}+6 y^{90}+3 y^{89}+6 y^{88}+4 y^{87}+3 y^{86}+5 y^{85}+4 y^{84}+7 y^{83}+3 y^{82}+4 y^{81}+3 y^{80}+y^{79}+y^{78}+y^{76}+y^{73} \quad,
$$
hence roughly a normal distribution around $100$.

This (and similar runs) gave us faith in our implementation.

{\bf Conclusion}

By using the problem of enumerating rectangular tilings of grid graphs with a fixed width, pioneered by David Klarner and Spyros Magliveras [KS], and
continued brilliantly by Joshua Smith and Helena Verrill [SV], as a {\it case study}, we demonstrated the power
of experimental mathematics and symbolic computation.

{\bf Encore}

Unrelated to the mathematics, but still fun, is the Maple package

{\tt \url{https://sites.math.rutgers.edu/~zeilberg/tokhniot/Recto.txt}} \quad,

that creates and solves Recto puzzles.

{\bf References}

[KM] David A. Klarner and Spyros S. Magliveras, {\it  The number of tilings of a block with blocks}, European Journal of Combinatorics {\bf 9} (1988), 317-330. \hfill\break
{\tt https://www.sciencedirect.com/science/article/pii/S0195669888800623}

[SV] Joshua Smith and Helena Verrill, {\it On dividing a rectangle into rectangles},\hfill\break
{\tt https://oeis.org/A116694/a116694.pdf} \quad.

[Sl] Neil Sloane, {\it The On-Line encyclopedia of integer sequences}, {\tt https://oeis.org/} . 

[W] Herbert S. Wilf, A unified setting for sequencing, ranking, and selection algorithms for combinatorial objects, Advances in Math 24 (1977) , 281-291. \hfill\break
{\tt https://www2.math.upenn.edu/~wilf/website/Unified\%20setting.pdf} \quad .

[Z] Doron Zeilberger, {\it The C-finite ansatz}, Ramanujan journal {\bf 31} (2013), 23-32. \hfill\break
{\tt https://sites.math.rutgers.edu/~zeilberg/mamarim/mamarimhtml/cfinite.html} \quad .

\bigskip
\hrule
\bigskip
Pablo Blanco, Natalya Ter-Saakov and Doron Zeilberger, Department of Mathematics, Rutgers University (New Brunswick), Hill Center-Busch Campus, 110 Frelinghuysen
Rd., Piscataway, NJ 08854-8019, USA. \hfill\break
Emails: {\tt  pablancoh at aol dot com} , {\tt  nt399 at rutgers dot edu} ,{\tt DoronZeil at gmail  dot com}   \quad .
\bigskip
Robert Dougherty-Bliss, Department of Mathematics, Dartmouth College, {\tt robert dot w dot bliss at gmail dot com} 
\bigskip

{\bf June 1, 2026} 

\end

\end{plain}
\end{document}